\documentclass[a4paper,12pt]{article}

\usepackage[usenames,dvipsnames]{pstricks}

\usepackage{amscd,amssymb,amsmath,graphicx,verbatim,hyperref}
\usepackage[OT1,T1]{fontenc}

\usepackage[all, cmtip]{xy}
\usepackage{pdfpages}

\usepackage{amsmath,amsthm,amsfonts,mathabx,tikz,tikz-cd,stmaryrd,url,mathtools, mathrsfs, tensor}
\usepackage{hyperref}
\hypersetup{
	colorlinks,
	citecolor=blue,
	filecolor=black,
	linkcolor=blue,
	urlcolor=red
}

\usetikzlibrary{matrix,arrows,decorations.pathmorphing}

\newcounter{Definitioncount}
\newtheorem{theorem}{Theorem}
\newtheorem{lemma}[theorem]{Lemma}
\newtheorem{proposition}[theorem]{Proposition}
\newtheorem{corollary}[theorem]{Corollary}

\theoremstyle{definition}
\newtheorem{definition}{Definition}
\newtheorem{remark}[theorem]{Remark}
\newtheorem{example}{Example}

\newtheoremstyle{fact}{\bigskipamount}{\medskipamount}{\upshape}{}{\itshape}{. }{ }{Fact}
\theoremstyle{fact}

\newtheoremstyle{quest}{\bigskipamount}{\medskipamount}{\upshape}{}{\itshape}{. }{ }{Question}
\theoremstyle{quest}

\newtheoremstyle{step}{2\bigskipamount}{\medskipamount}{\upshape}{}{\itshape}{. }{ }{\underline{Step~\thestep}}
\theoremstyle{step}

\renewcommand{\thestep}{\arabic{step}}

\newcommand{\lra}{\longrightarrow}
\newcommand{\Ra}{\Rightarrow}

\newcommand{\ldual}[1]{\mathord{{\let\nolimits\relax\sideset{^\wedge}{}{#1}}}}
\newcommand{\laction}[2]{\mathord{{\let\nolimits\relax\sideset{^{#1}}{}{#2}}}}
\newcommand{\conj}[2]{\mathord{{\let\nolimits\relax\sideset{^{#1}}{}{#2}}}}

\newcommand{\ox}{\otimes}
\newcommand{\xra}{\xrightarrow}
\newcommand{\xla}{\xleftarrow}

\newcommand{\nra}{\nrightarrow}

\usepackage[utf8]{inputenc}
\DeclareFontFamily{U}{min}{}
\DeclareFontShape{U}{min}{m}{n}{<-> udmj30}{}
\newcommand\yo{\!\text{\usefont{U}{min}{m}{n}\symbol{'207}}\!}

\def\CA{{\mathscr A}}
\def\CB{{\mathscr B}}
\def\CC{{\mathscr C}}
\def\CD{{\mathscr D}}

\def\CI{{\mathscr I}}

\def\CK{{\mathscr K}}

\def\CV{{\mathscr V}}

\def\CX{{\mathscr X}}

\makeatletter
\newcommand*\bigcdot{\mathpalette\bigcdot@{.5}}
\newcommand*\bigcdot@[2]{\mathbin{\vcenter{\hbox{\scalebox{#2}{$\m@th#1\bullet$}}}}}
\makeatother

\newcommand{\twocong}[2][0.5]{\ar@{}[#2] \save ?(#1)*{\cong}\restore}
\newcommand{\twoeq}[2][0.5]{\ar@{}[#2] \save ?(#1)*{=}\restore}
\newcommand{\ltwocell}[3][0.5]{\ar@{}[#2] \ar@{=>}?(#1)+/r 0.2cm/;?(#1)+/l 0.2cm/^{#3}}
\newcommand{\rtwocell}[3][0.5]{\ar@{}[#2] \ar@{=>}?(#1)+/l 0.2cm/;?(#1)+/r 0.2cm/^{#3}}
\newcommand{\utwocell}[3][0.5]{\ar@{}[#2] \ar@{=>}?(#1)+/d  0.2cm/;?(#1)+/u 0.2cm/_{#3}}
\newcommand{\dtwocell}[3][0.5]{\ar@{}[#2] \ar@{=>}?(#1)+/u  0.2cm/;?(#1)+/d 0.2cm/^{#3}}
\newcommand{\ultwocell}[3][0.5]{\ar@{}[#2] \ar@{=>}?(#1)+/dr  0.2cm/;?(#1)+/ul 0.2cm/^{#3}}
\newcommand{\urtwocell}[3][0.5]{\ar@{}[#2] \ar@{=>}?(#1)+/dl  0.2cm/;?(#1)+/ur 0.2cm/^{#3}}
\newcommand{\dltwocell}[3][0.5]{\ar@{}[#2] \ar@{=>}?(#1)+/ur  0.2cm/;?(#1)+/dl 0.2cm/^{#3}}
\newcommand{\drtwocell}[3][0.5]{\ar@{}[#2] \ar@{=>}?(#1)+/ul  0.2cm/;?(#1)+/dr 0.2cm/^{#3}}

\begin{document}
\author{Ross Street \footnote{The author gratefully acknowledges the support of Australian Research Council Discovery Grant DP190102432.} \\ 
\small{Centre of Australian Category Theory} \\
\small{Macquarie University, NSW 2109 Australia} \\
\small{<ross.street@mq.edu.au>}
}
\title{Wood fusion for magmal comonads}

\date{\today}
\maketitle


\begin{abstract}
The goal is to show how a 1978 paper of Richard Wood on monoidal comonads and exponentiation relates to more recent publications such as Pastro et alia \cite{97} and Brugui{\'e}res et alia \cite{Bruguieres2011}.
In the process, we mildly extend the ideas to procomonads in a magmal setting and suggest it also works for algebras for any club in the sense of Max Kelly \cite{Kelly1972(a), Kelly1972(b)}.
\end{abstract}
\noindent {\small{\emph{2010 Mathematics Subject Classification:} 16E45; 16D90; 18A40}
\tableofcontents

\section{Introduction}
This paper began when I moved offices and uncovering Richard Wood's paper \cite{Wood1978}
in which he considers a closed monad $G$ on a closed category $\CC$. He
provides necessary and sufficient conditions for the natural promonoidal structure on the
category $\CC^G$ of Eilenberg-Moore coalgebras for $G$ to be closed (that is, to be
representable by an internal hom). He also provides necessary and sufficient conditions 
for $\CC^G$ to be closed in that way and for the underlying functor $\mathrm{und} : \CC^G\to \CC$
to be strong closed (that is, to preserve the internal hom).   

The process involves morphisms
\begin{eqnarray}\label{WoodFusionIntro}
w_{\upsilon, Z} : = G[Y,Z] \xra{G^{\ell}_2} [GY,GZ] \xra{[\upsilon,1]} [Y,GZ] 
\end{eqnarray}  
for $Y\xra{\upsilon}GY$ an object of $\CC^G$ and $Z$ an object of $\CC$.
Here I am denoting the closed structure on $G$ by $G^{\ell}_2$ where the closed
structure on $\CC$ is what I think of as the ``left'' kind. These are what I will call
{\em Wood fusion morphisms}.  

An aim of the present paper is to relate Wood fusion to the fusion occurring in more recent papers
such as \cite{Street1998, Bruguieres2011, 103, Pastro2012}. 

Like Wood, I will write for monoidal closed $\CC$. Actually, it became clear 
that the story applies to more general categorical structures. In particular,
we could work with algebras $\CC$ for a club in the sense of Kelly \cite{Kelly1972(a), Kelly1972(b)}.
It is about the functors providing the operations (such as binary tensor product), while
the structural natural transformations (such as associativity constraints) and axioms thereon 
(such as the Mac Lane pentagon) carry over automatically to the constructions.
However, rather than review the notion of club and to keep the notation simple, I have written
for the case of {\em magmal categories}; that is, categories $\CC$ equipped with a binary functorial
operation $\CC\times \CC \to \CC$ which we still call and write as a {\em tensor product}.
Indeed, as is our custom, we will be a little more general and work with hom $\CV$-enriched categories throughout.

The unification of the monad and comonad cases is achieved by working with a magmal {\em procomonad}
$\Gamma$ on a magmal $\CC$.

Appendices on adjoint lifting theorems of \cite{Appel1965, DubucAT, Johnstone1975, Keigher1975} are added to provide a perspective suggested by Richard Garner and Stephen Lack.
 
\section{Magmal comonads and cloaks}

We work with categories enriched in a suitably complete and cocomplete, symmetric, monoidal category $\CV$ as base for enrichment (see Kelly \cite{KellyBook}). The tensor of $\CV$ is denoted by $U\ox V$
and the tensor unit by $I$. At times, we omit the prefix ``$\CV$-'', taking it for granted.

\begin{definition}\label{defmag} 
A $\CV$-category $\CC$ is {\em magmal} when it is equipped with a $\CV$-functor
$\ox : \CC\ox \CC \to \CC$ (this overworked notation, written as usual between the arguments). A $\CV$-functor $S : \CC\to \CD$
between magmal $\CV$-categories is called {\em magmal} when it is equipped with a $\CV$-natural
transformation $S_2$ as in \eqref{magfun}. 
\begin{equation}\label{magfun}
 \begin{aligned}
  \xymatrix{
\CC\ox \CC \ar[d]_{S\ox S}^(0.5){\phantom{aaaaaa}}="1" \ar[rr]^{\ox}  && \CC \ar[d]^{S}_(0.5){\phantom{aaaaaa}}="2" \ar@{=>}"1";"2"^-{S_2}
\\
\CD\ox \CD \ar[rr]_-{\ox} && \CD 
}
 \end{aligned}
\end{equation}
A $\CV$-natural transformation $\theta : S\Ra T : \CC\to \CD$
between magmal $\CV$-functors is called {\em magmal} when it satisfies the equation
$\theta_{C\ox C'} \circ S_{2; C, C'} = T_{2; SC,SC'}\circ (\theta_C\ox \theta_{C'})$
for all $C, C'\in \CC$. A comonad $G=(G,\varepsilon, \delta)$ is {\em magmal} when
$G$, $\varepsilon : G\Ra 1_{\CC}$ and $\delta : G\Ra G\circ G$ are all magmal.   
\end{definition}

\begin{definition}\label{defncloaked}
Suppose $\CC$ is a magmal $\CV$-category and $Y, Z \in \CC$.
We say {\em $Z$ is left cloaked by $Y$} when there is an object $[Y,Z]\in \CC$ and a morphism $$\mathrm{ev}^Y_Z : [Y,Z]\ox Y\to Z $$ such that the family
$$\CC(X,[Y,Z]) \lra \CC(X\ox Y,Z) \ ,$$
$\CV$-natural in $X$ obtained by the Yoneda Lemma, consists of isomorphisms.
In the language of Appendix~\ref{liftings}, $[Y,Z] = \mathrm{rif}(-\ox Y,Z)$ and $\mathrm{ev}^Y_Z = \varepsilon_Z^{-\ox Y}$ in the 2-category $\CV\text{-}\mathrm{Cat}$.  
We say $[Y,Z]$ is the {\em left cloak of $Z$ by $Y$};
if these exist for all $Z$, we have the components
\begin{eqnarray*}
\mathrm{ev}^Y_Z : [Y,Z]\ox Y\to Z \ \text{ and } \ \mathrm{ve}_X^Y : X\to [Y,X\ox Y]
\end{eqnarray*}
of the counit and unit of an adjunction $-\ox Y\dashv [Y,-]$.
If left cloaks $[Y,Z]$ exist for all $Y, Z\in \CC$, we say {\em $\CC$ is left cloakal}.      
\end{definition}

\begin{proposition}\label{magmalcloakal}
Suppose $S : \CC\to \CD$ is a $\CV$-functor between magmal $\CV$-categories and $Y\in \CC$.
Suppose $\CC$ admits left cloaks by $Y$ and $\CD$ admits left cloaks by $SY$.
There is a bijection between families
$$S_{2; X,Y} : SX\ox SY \to S(X\ox Y)$$
$\CV$-natural in $X$ and families
$$S_{2; Y,Z}^{\ell} : S[Y,Z]\to [SY,SZ]$$
$\CV$-natural in $Z$. 
\end{proposition}
\begin{proof}
This is an application of the mate bijection (in the sense of \cite{KelSt1974}) between $\CV$-natural $(-\ox SY)\circ S\xra{S_{2; -,Y}} S\circ (-\ox Y)$ and
$S\circ [Y,-]\xra{S_{2; Y,-}^{\ell}} [SY,-]\circ S$ under the adjunctions $-\ox Y\dashv [Y,-]$
and $-\ox SY\dashv [SY,-]$. 
\end{proof}

Obviously every monoidal category, functor, and natural transformation has an underlying magmal structure.  

Let $G$ be a magmal comonad on a magmal $\CV$-category $\CC$.
So we have $\CV$-natural families  $\varepsilon_X : GX\to X$,
$\delta_{X} : GX\to GGX$ and $G_{2; X,Y} : GX\ox GY\to G(X\ox Y)$. 

As for monoidal $\CV$-categories (see \cite{Wood1978, Moerd02, McC02}), the category $\CC^G$ of Eilenberg-Moore ${G}$-coalgebras becomes magmal with tensor product
\begin{eqnarray}\label{magcoalg}
(X\xra{\xi}GX)\ox (Y\xra{\upsilon}GY) = (X\ox Y\xra{\xi\ox \upsilon} GX\ox GY\xra{G_{2; X,Y}}G(X\ox Y)) \ .
\end{eqnarray}
This gives the construction of coalgebras for comonads (in the sense of \cite{3}) in the 2-category $\mathrm{Mag}(\CV\text{-}\mathrm{Cat})$ of magmal $\CV$-categories, magmal $\CV$-functors, and
magmal $\CV$-natural transformations. 

The following result is an application of Proposition~\ref{rifwithM} and yet we still give a direct proof within the present context.
\begin{lemma}\label{magequalizer}
Let $G$ be a magmal comonad on a magmal $\CV$-category $\CC$.
Suppose cloaks $[(Y,\upsilon),(GZ,\delta_Z)]$ and $[(Y,\upsilon),(G^2Z,\delta_{GZ})]$ exist in $\CC^G$
for the objects $(Y,\upsilon)$ and $(Z,\zeta)$ of $\CC^G$.
Then the cloak $[(Y,\upsilon),(Z,\zeta)]$ exists if and only if the parallel pair 
\begin{eqnarray}\label{011}
\xymatrix @R-3mm {
 [(Y,\upsilon),(GZ,\delta_Z)] \ar@<1.5ex>[rr]^{[1,\delta_Z]}  \ar@<-1.5ex>[rr]_{[1,G\zeta]}  && [(Y,\upsilon), (G^2Z,\delta_{GZ})] \ .}
\end{eqnarray}
has an equalizer $(E,\kappa)\xra{k}[(Y,\upsilon),(GZ,\delta_Z)]$ in $\CC^G$. 
If that holds, then $$[(Y,\upsilon),(Z,\zeta)]= (E,\kappa)$$ 
and the counit $\mathrm{ev}^{(Y,\upsilon)}_{(Z,\zeta)}$ is determined by commutativity of the square \eqref{evalcompat}.
\begin{eqnarray}\label{evalcompat}
\begin{aligned}
\xymatrix{
(E,\kappa)\ox (Y,\upsilon) \ar[rr]^-{\mathrm{ev}^{(Y,\upsilon)}_{(Z,\zeta)}} \ar[d]_-{k\ox 1} && (Z,\zeta) \ar[d]^-{\zeta} \\
[(Y,\upsilon),(GZ,\delta_Z)]\ox (Y,\upsilon) \ar[rr]_-{\mathrm{ev}^{(Y,\upsilon)}_{(GZ,\delta_Z)}} && (GZ,\delta_Z)}
\end{aligned}
\end{eqnarray}
\end{lemma}
\begin{proof} We have the equalizer
\begin{eqnarray}\label{110}
\xymatrix @R-3mm {
(Z,\zeta) \ar@<0.0ex>[rr]^{\zeta} && (GZ,\delta_Z) \ar@<1.5ex>[rr]^{\delta_Z}  \ar@<-1.5ex>[rr]_{G\zeta}  && (G^2Z,\delta_{GZ}) 
}\end{eqnarray}
in $\CC^G$. The ``only if'' direction is then clear since right adjointness of the functor $[(X,\xi),-]$ implies it preserves the equalizer \eqref{110}.
For the other direction, note that morphisms $f : (A,\alpha)\to (E,\kappa)$ are in bijection (via composition with $k$) with morphisms $g : (A,\alpha)\to [(Y,\upsilon),(GZ,\delta_Z)]$ such that $[1,\delta_Z]\circ g = [1,G\zeta]\circ g$ and so with morphisms $h : (A,\alpha)\ox (Y,\upsilon) \to (GZ,\delta_Z)$ such that $\delta_Z\circ h = G\zeta\circ h$. Using the equalizer \eqref{110}, we see that $(E,\kappa)$ is the stated cloak and that the last sentence of the lemma holds.      
\end{proof}

The following lemma is essentially in \cite{Wood1978}.

\begin{lemma}\label{cloak for cofree}
Let $G$ be a magmal comonad on a magmal $\CV$-category $\CC$.
Suppose the cloaks $[Y,Z]$ and $[GY,GZ]$ exist in $\CC$.
Then a cloak of the cofree object $(GZ,\delta_Z)$ by any object $(Y,\upsilon)$ in $\CC^{ {G}}$
is given by the cofree object $(G[Y,Z],\delta_{[Y,Z]})$ with the composite  
\begin{eqnarray*}
G[Y,Z]\ox Y \xra{G^{\ell}_2\ox \xi} [GY,GZ]\ox GY\xra{\mathrm{ev}^{GY}_{GZ}} GZ
\end{eqnarray*}
as $\mathrm{ev}^{(Y,\upsilon)}_{(GZ,\delta_Z)}$.
\end{lemma}
\begin{proof} We have the chain of isomorphisms
\begin{eqnarray*}
\lefteqn{\CC^G((A,\alpha),(G[Y,Z],\delta_{[Y,Z]})) } \\
& \cong & \CC(A,[Y,Z]) \\
& \cong & \CC(A\ox Y,Z) \\
& \cong & \CC^G((A\ox Y, \phi (\alpha\ox \upsilon)),(GZ,\delta_Z)) \\
& \cong & \CC^G((A,\alpha)\ox (Y,\upsilon),(GZ,\delta_Z))
\end{eqnarray*}
which by Yoneda is induced by the evaluation morphism given by the clockwise path
around the diagram \eqref{equalevals}.
\begin{eqnarray}\label{equalevals}
\begin{aligned}
\xymatrix{
G[Y,Z]\ox Y \ar[rrdd]_-{G^{\ell}_2\ox \upsilon} \ar[rr]^-{\delta_{[Y,Z]}\ox \upsilon} \ar[rrd]^-{1\ox \upsilon} && GG[Y,Z]\ox GY \ar[rr]^-{G_2} \ar[d]^-{G\varepsilon_{[Y,Z]}\ox 1} && G(G[Y,Z]\ox Y) \ar[d]^{G(\varepsilon_{[Y,Z]}\ox 1)} \\
 && G[Y,Z]\ox GY \ar[rr]_-{G_2} \ar[d]^-{G^{\ell}_2 \ox 1} && G([Y,Z]\ox Y) \ar[d]^-{G\mathrm{ev}^Y_Z} \\
&& [GY,GZ]\ox GY \ar[rr]_-{\mathrm{ev}^{GY}_{GZ}} && GZ}
\end{aligned}
\end{eqnarray}
Diagram \eqref{equalevals} thus shows $\mathrm{ev}^{(Y,\upsilon)}_{(GZ,\delta_Z)}$ to be the composite stated in the lemma.
\end{proof}

\section{Wood fusion morphisms}

Notice that, if the cloak $[Y,GZ]$ exists, then the evaluation in Lemma \ref{cloak for cofree} corresponds,
under the universal property of cloaks, to the composite
\begin{eqnarray}\label{WoodFusion}
w_{\upsilon, Z} : = G[Y,Z] \xra{G^{\ell}_2} [GY,GZ] \xra{[\upsilon,1]} [Y,GZ] \ .
\end{eqnarray} 
 \begin{definition}\label{defwoodfusion}
 {\em The $w_{\upsilon, Z}$ of \eqref{WoodFusion} are called {\em Wood $G$-fusion morphisms}.} 
  \end{definition}
  Notice that the Wood fusion morphism for cofree coalgebras occurs in the construction of a new skew-closed structure using a closed comonad; see Proposition 3 of \cite{116}.
  
  We will be interested in when the Wood fusion morphisms are invertible.
  As with ordinary fusion (recalled in Section~\ref{Ffom}), invertibility for an arbitrary 
  $G$-coalgebra follows from invertibility for cofree $G$-coalgebras. 
  
  \begin{proposition}\label{cofreesuffice} For a $G$-coalgebra $(Y,\upsilon)$ and any
  object $Z$, the Wood fusion morphism $w_{\upsilon, Z}$ is invertible if $w_{\delta_Y, Z}$
  is invertible and $w_{\delta_{GY}, Z}$ is an epimorphism.   
\end{proposition}
\begin{proof}
For any coalgebra $(Z,\zeta)$, we have an equalizer of the form
\begin{eqnarray}\label{122}
\xymatrix @R-3mm {
Z \ar@<0.0ex>[rr]^{\zeta} && GZ \ar@<1.5ex>[rr]^{\delta_Z}  \ar@<-1.5ex>[rr]_{G\zeta}  && G^2Z  
}\end{eqnarray} 
which is preserved by all functors (that is, it is an absolute equalizer).
It follows that the rows of \eqref{3coeq} are coequalizers.
The vertical morphisms give two composable morphisms of coequalizer diagrams
(that is, the appropriate diagrams commute using naturality of $G_2$, and coassociativity properties of $\upsilon$ and $\delta$).  
\begin{eqnarray}\label{3coeq}
\begin{aligned}
\xymatrix @R-3mm {
G[G^2Y,Z]\ar[d]_-{G^{\ell}_2} \ar@<1.5ex>[rr]^{G[\delta_Y,1]}  \ar@<-1.5ex>[rr]^{G[G\upsilon,1]} && G[GY,Z]\ar[d]^-{G^{\ell}_2} \ar@<0.0ex>[rr]^{G[\upsilon,1]} && G[Y,Z] \ar[d]^-{G^{\ell}_2}\\ 
[G^3Y,GZ]\ar[d]_-{[\delta_{GY},1]} \ar@<1.5ex>[rr]^{[G\delta_Y,1]}  \ar@<-1.5ex>[rr]^{[G^2\upsilon,1]} && [G^2Y,GZ]\ar[d]^-{[\delta_{Y},1]} \ar@<0.0ex>[rr]^{[G\upsilon,1]} && [GY,GZ] \ar[d]^-{[\upsilon,1]}\\
[G^2Y,GZ] \ar@<1.5ex>[rr]^{[\delta_Y,1]}  \ar@<-1.5ex>[rr]^{[G\upsilon,1]} && [GY,GZ] \ar@<0.0ex>[rr]^{[\upsilon,1]} && [Y,GZ]}
\end{aligned}
\end{eqnarray}
If the middle vertical composite is invertible and the left vertical composite is an
epimorphism then the right vertical composite is invertible.    
\end{proof}

\begin{definition}
 {\em The closed comonad $G$ on $\CC$ is {\em Hopf-Wood} when the Wood fusion morphisms $w_{\upsilon, Z}$ are all invertible. By Proposition~\ref{cofreesuffice}, it
 suffices to know that all $w_{\delta_Y, Z}$ are invertible; then the property can
 be expressed without reference to Eilenberg-Moore coalgebras.}  
  \end{definition}
  
  \begin{lemma}\label{transportedpair} Under the conditions of Lemma~\ref{cloak for cofree}, the parallel
  pair \eqref{011} is isomorphic to the parallel pair
  \begin{eqnarray}\label{021}
\xymatrix @R-3mm {
 (G[Y,Z],\delta_{[Y,Z]}) \ar@<1.5ex>[rr]^{Gw_{\upsilon, Z}\circ \delta_{[Y,Z]}}  \ar@<-1.5ex>[rr]_{G[1,\zeta]}  && (G[Y,GZ],\delta_{[Y,GZ]}) \ .}
\end{eqnarray}
\end{lemma}
  
\begin{definition} A strong magmal functor $K : \CA\to \CC$ is said to {\em create the cloak of $B$ by $A$} in $\CA$ when
the cloak of $KB$ by $KA$ exists in $\CC$, and there exist $H$ in $\CA$ and $\tau : KH\cong [KA,KB]$ in $\CC$ 
with the following two properties:
\begin{itemize}
\item[(i)] there exists a unique morphism $\bar{e} : H\ox A \to B$ such that the square
\begin{equation}\label{cloakcreation}
\begin{aligned}
\xymatrix{
KH\ox KA \ar[rr]^-{\tau\ox 1} \ar[d]_-{\phi} && [KA,KB]\ox KA \ar[d]^-{e_{KB}^{KA}} \\
K(H\ox A) \ar[rr]_-{K\bar{e}} && KB}
\end{aligned}
\end{equation}
commutes;
\item[(ii)] the object $H$ with $\bar{e}$ is a cloak for $B$ by $A$.   
\end{itemize} 
\end{definition}
The following lemma is straightforward.
\begin{lemma}\label{pbcloakcreation} Suppose \eqref{pbstrongmag} is a pullback of magmal $\CV$-categories
and strong magmal $\CV$-functors with $W$ fully faithful. Suppose $A, B\in \CA$ are such that $K'$ creates the
cloak of $VB$ by $VA$. If $[VA,VB] \cong VH$ for some $H$ then $K$ creates the cloak of $B$ by $A$.
 \begin{equation}\label{pbstrongmag}
\begin{aligned}
\xymatrix{
\CA \ar[rr]^-{V} \ar[d]_-{K} && \CA' \ar[d]^-{K'} \\
\CC \ar[rr]_-{W} && \CC'}
\end{aligned}
\end{equation}
\end{lemma}
\begin{lemma}\label{restricted creation} Suppose $[Y,Z]$ and $[GY,GZ]$ exist in $\CC$.
The forgetful functor $$\mathrm{und}_G : \CC^{ {G}} \to \CC$$ creates the exponential $[(Y,\upsilon),(GZ,\delta_Z)]$ of
Lemma \ref{cloak for cofree} if and only if $[Y,GZ]$ exists in $\CC$ and the Wood $G$-fusion morphism \eqref{WoodFusion} is invertible.
\end{lemma}
\begin{proof} Assume $[Y,GZ]$ exists in $\CC$ and the Wood $G$-fusion morphism \eqref{WoodFusion} is invertible.
In the definition of creation, take $H = (G[Y,Z],\delta_{[Y,Z]})$ and $\tau = w_{\upsilon, Z}$. 
Commutativity of diagram \eqref{cloakcreation} means, in this case, that $\bar{e}$ must be the clockwise route around the diagram
\begin{eqnarray*}
\xymatrix{
G[Y,Z]\ox Y \ar[rr]^-{G^{\ell}_2\ox 1} \ar[d]_-{1\ox\upsilon} && [GY,GZ]\ox Y \ar[rr]^-{[\upsilon,1]\ox 1} \ar[d]^-{1\ox\upsilon} && [Y,GZ]\ox Y \ar[d]^-{\mathrm{ev}^Y_{GZ}} \\
G[Y,Z]\ox GY \ar[rr]_-{G^{\ell}_2\ox 1}  && [GY,GZ]\ox GY \ar[rr]_-{\mathrm{ev}^{GY}_{GZ}}  && GZ}
\end{eqnarray*}
and the diagram shows that $\bar{e}$ is the evaluation displayed in Lemma~\ref{cloak for cofree}.
By Lemma~\ref{cloak for cofree} we have what we need for ``if''.

Now suppose $\mathrm{und}_G : \CC^{ {G}} \to \CC$ creates the exponential $[(Y,\upsilon),(GZ,\delta_Z)]$ which we know from 
Lemma~\ref{cloak for cofree} is the object $(G[Y,Z], \delta_{[Y,Z]})$ with the evaluation $\bar{e}$ as displayed in that lemma.
Since $\mathrm{und}_G(GZ,\delta_Z)= GZ$ and $\mathrm{und}_G(Y,\upsilon)=Y$, we have the existence of $[Y,GZ]$
and that there is an isomorphism $\tau : G[Y,Z]\cong [Y,GZ]$ such that $\bar{e}=\mathrm{ev}^X_{GY}\circ (\tau\ox 1)$. 
This last equation means that $\tau$ corresponds to $\bar{e}$ under the universal property of $[Y,GZ]$;
that is, $\tau = w_{\upsilon, Z}$, which proves ``only if''.    
\end{proof}

\begin{proposition}\label{magcomoncloaks} Let $ {G}$ be a magmal comonad on a magmal $\CV$-category $\CC$.
Suppose $(Y,\upsilon)\in \CC^{ {G}}$ is such that $[Y,Z]$ and $[GY,GZ]$ exist for all $Z\in \CC$. 
The Wood $G$-fusion morphisms $w_{\upsilon, Z}$ are invertible for
all $Z$ if and only if $\mathrm{und}_G : \CC^{ {G}} \to \CC$ creates all cloaks by $(Y,\upsilon)$.
In this case, for any $(Z,\zeta)\in \CC^{ {G}}$, 
$$[(Y,\upsilon),(Z,\zeta)] \cong ([Y,Z], w_{\upsilon,Z}^{-1}\circ [1,\zeta]) \ .$$ 
\end{proposition}
\begin{proof} Lemma~\ref{restricted creation} gives ``if''. Suppose all $w_{(Y,\upsilon), Y}$ are invertible.
We will use facts involved in the Beck monadicity theorem \cite{CWM} in dual form for comonads.
We have the cosplit equalizer \eqref{122}.
From Lemma~\ref{restricted creation}, the parallel pair in equalizer \eqref{011}
is taken by $U_G$ to the cosplit pair
 \begin{eqnarray}\label{114}
 \xymatrix @R-3mm {
[Y,GZ] \ar@<1.5ex>[rr]^{[1,\delta_Z]}  \ar@<-1.5ex>[rr]_{[1,G\zeta]}  && [Y, GGZ] }
\end{eqnarray}
which, by applying $[Y,-]$ to \eqref{122}, has the cosplit equalizer
\begin{eqnarray}\label{115}
\xymatrix @R-3mm {
[Y,Z] \ar@<0.0ex>[rr]^{[1,\zeta]} && [Y,GZ] \ar@<1.5ex>[rr]^{[1,\delta_Z]}  \ar@<-1.5ex>[rr]_{[1,G\zeta]}  && [Y, GGZ] \ . 
}\end{eqnarray}
Since $\mathrm{und}_G$ is comonadic, there exists a unique $G$-coalgebra $([Y,Z],\kappa)$ and an equalizer
\begin{eqnarray}\label{114}
\xymatrix @R-3mm {
 ([Y,Z],\kappa) \ar@<0.0ex>[rr]^{[1,\zeta]} && ([Y,GZ],\kappa_1) \ar@<1.5ex>[rr]^{[1,\delta_Z]}  \ar@<-1.5ex>[rr]_{[1,G\zeta]}  && ([Y, G^2Z],\kappa_2)}
\end{eqnarray}
in $\CC^G$, where (using Lemma~\ref{transportedpair}) the coactions $\kappa_1$ and $\kappa_2$ are transported from the coactions
$\delta_{[Y,Z]}$ on $G[Y,Z]$ and $\delta_{[Y,GZ]}$ on $G[Y,GZ]$ under the invertible Wood fusion morphisms. So we have condition (ii) for $\mathrm{und}_G$ to create the cloak. For condition (i),
note that commutativity of \eqref{evalcompat} in Lemma~\ref{magequalizer} with $k = [1,\zeta]$ shows
that $\mathrm{ev}^Y_Z : ([Y,Z],\kappa)\ox (Y,\upsilon)\to (Z,\zeta)$ is the G-coalgebra morphism for the unique solution to diagram \eqref{cloakcreation}. 
\end{proof}

\section{Fusion for opmagmal monads}\label{Ffom}

Let $T$ be an opmagmal monad on the magmal $\CV$-category $\CC$.
The monad structure involves a unit $\eta : 1_{\CC}\to T$ and a multiplication
$\mu : TT\to T$. The opmagmal structure involves a natural family of morphisms
$T_{2;X,Y} : T(X\ox Y)\to TX\ox TY$. 
We denote the magmal category of Eilenberg-Moore $T$-algebras by $\CC^T$ 
with strong magmal forgetful functor $\mathrm{und}_T : \CC^T\to \CC$. 
The tensor product for $\CC^T$ is defined by  
$$(TX\xra{\alpha}X)\ox (TY\xra{\beta}Y) = \left(T(X\ox Y)\xra{T_2} TX\ox TY\xra{\alpha\ox\beta}X\ox Y\right) \ .$$

For $X\in \CC$ and $(Y,\beta)\in \CC^T$, we call the composite
$v = v_{X,\beta}$:
\begin{eqnarray}\label{Tfusionmorph}
 T(X\ox Y)\xra{T_2}  TX\ox TY\xra{1\ox \beta} TX\ox Y
\end{eqnarray}
a {\em $T$-fusion morphism} (as featured in \cite{Bruguieres2011}).

Suppose $T : \CC\to \CC$ has a right adjoint functor $G$. As discussed in \cite{EM1965},
$G$ becomes a comonad on $\CC$ and there is an isomorphism of categories $\CC^T\cong \CC^G$
over $\CC$. These matters involve the calculus of mates (in the sense of \cite{KelSt1974}) as does
the fact that $G$ becomes a monoidal comonad and the isomorphism $\CC^T\cong \CC^G$ becomes strong monoidal.  

\begin{proposition}\label{adjointprop} Suppose $T$ is an opmagmal monad on the magmal $\CV$-category $\CC$.
Suppose $G$ is a right adjoint magmal comonad for $T$. Let $(Y,\beta)\in \CC^T$ correspond
to $(Y,\upsilon)\in \CC^G$. The $T$-fusion morphism $v_{X,\beta}$ is invertible for all $X\in \CC$ if
and only if the Wood $G$-fusion morphism $w_{\upsilon,Z}$ is invertible for all $Z\in \CC$. 
\end{proposition}
\begin{proof}
Apply the Yoneda Lemma to the following commutative diagram where $\sigma : TG\to 1_{\CC}$ is
the counit of $T\dashv G$.
\begin{eqnarray*}
\xymatrix{
\CC(TX\ox Y,Z) \ar[rr]^-{\cong} \ar[dd]_-{\CC(1\ox\beta ,1)} \ar[rd]_-{\CC(1\ox\sigma_Y,1)} && \CC(X, G[Y,Z]) \ar[dd]^-{\CC(1,G_2^{\ell})}\\
& \CC(TX\ox TGY,Z) \ar[dd]^-{\CC(G_2 , 1)} \ar[ld]^-{\CC(1\ox T\upsilon ,1)}  & \\ 
\CC(TX\ox TY,Z) \ar[dd]_-{\CC(G_2 , 1)} && \CC(X,[GY,GZ]) \ar[dd]^-{\CC(1,[\upsilon,1])} \\
& \CC(T(X\ox GY),Z) \ar[ru]_-{\cong} \ar[ld]_-{\CC(T(1\ox\upsilon),1)} & \\
\CC(T(X\ox Y),Z) \ar[rr]_-{\cong} & & \CC(X,[Y,GZ])}
\end{eqnarray*}
where $\sigma : TG\to 1_{\CC}$ is the counit of $T\dashv G$.
\end{proof}

\begin{example} Let $\CC$ be a braided closed monoidal category. Let $H$ be a monoid in 
the monoidal category of comagma in $\CC$.
Then $-\ox H : \CC\to \CC$ is an opmagmal monad with right adjoint 
$[H,-] : \CC\to \CC$. 
Proposition~\ref{adjointprop} relates Wood fusion for the magmal comonad 
$[H,-]$ with the fusion morphism 
$$H\ox H \xra{\delta\ox 1} H\ox H\ox H \xra{1\ox \mu} H\ox H$$
for $H$.  We say $H$ is Hopf when its fusion morphism is invertible. 
This is equivalent to $-\ox H$ Hopf and to $[H,-]$ Hopf-Wood.  
\end{example}

\section{Procomonads}\label{Pc}

Let $\CV$ be a symmetric closed monoidal category which is complete and cocomplete.
Let $\mathfrak{M} = \CV\text{-}\mathrm{Mod}$ be the bicategory of $\CV$-categories and $\CV$-modules
in the terminology of \cite{88, 60} and elsewhere; modules are also called ``bimodules'' by 
Lawvere \cite{LawMetricSp}, and first ``profunctors'' \cite{Ben1967} and then 
``distributors'' \cite{Ben1973} by B\'enabou. 
The bicategory $\mathfrak{M}$ has homs enriched in $\CV\text{-}\mathrm{Cat}$;
we equate the hom $\mathfrak{M}(\CA,\CB)$ with the $\CV$-functor $\CV$-category $[\CB^{\mathrm{op}}\ox \CA,\CV]$.
Composition of $\CV$-modules $M : \CA\nrightarrow \CB$ and $N : \CB\nrightarrow \CC$ is defined by
coends
\begin{eqnarray*}
(N\circ M)(C,A) = \int^BM(B,A)\ox N(C,B) \ .
\end{eqnarray*}
Each $\CV$-functor $F : \CA\to \CB$ gives $\CV$-modules $F_*:\CA\nra \CB$ and 
$F^* : \CB\nra \CA$ with $F_*\dashv F^*$ in $\mathfrak{M}$; indeed, $F_*(B,A) = \CB(B,FA)$
and $F^*(A,B)=\CB(FA,B)$. A module $M :\CA\nra \CB$ is called {\em Cauchy} when it has a right
adjoint in $\mathfrak{M}$. A module $M :\CA\nra \CB$ is called {\em convergent} or {\em representable} when 
$M \cong F_*$ for some $\CV$-functor $F : \CA\to \CB$. 

We write $\CI$ for the $\CV$-category with one object $0$ and hom $\CI(0,0) = I$ (the tensor unit of $\CV$). Then $\mathfrak{M}(\CI,\CC)\cong [\CC^{\mathrm{op}},\CV]$,
the category of $\CV$-presheaves on $\CC$. Composition with $N\in \mathfrak{M}(\CB,\CC)$ transports to a left adjoint $\CV$-functor 
$$\widebar{N} : [\CB^{\mathrm{op}},\CV] \lra [\CC^{\mathrm{op}},\CV]$$ 
where 
\begin{eqnarray}\label{barring}
\widebar{N}(F)(C) = \int^BF(B)\ox N(C,B) \ \text{ so that } \ \widebar{N}(\CB(-,B)) \cong N(-, B) \ .
\end{eqnarray}
In fact, $N\mapsto \widebar{N}$ is the object function of a biequivalence between the 
bicategory $\mathfrak{M}$ and the 2-category $\mathfrak{P}$ of $\CV$-presheaf categories, left adjoint $\CV$-functors, and $\CV$-natural transformations.  

A {\em $\CV$-procomonad} is a comonad in $\mathfrak{M}$ and so consists of a $\CV$-category $\CC$, a $\CV$-module $\Gamma : \CC\nrightarrow\CC$, 
a $\CV$-natural transformation $\varepsilon : \Gamma \Ra 1_{\CC}$, and a $\CV$-natural transformation 
$\delta : \Gamma \Ra \Gamma\circ \Gamma$ satisfying the coassociativity and counital conditions.
We say $\Gamma = (\Gamma , \varepsilon, \delta)$ is a {\em $\CV$-procomonad on $\CC$}.  
For any $\CV$-category $\CA$, we obtain a $\CV$-comonad $\mathfrak{M}(1_{\CA}, \Gamma)$ on the $\CV$-category
$\mathfrak{M}(\CA,\CC)$. Then we have the $\CV$-category $\mathfrak{M}(\CA,\CC)^{\mathfrak{M}(1_{\CA}, \Gamma)}$ of Eilenberg-Moore
$\mathfrak{M}(1_{\CA}, \Gamma)$-coalgebras.
In particular, when $\CA = \CI$, we obtain a $\CV$-comonad $\widebar{\Gamma}$
on the presheaf $\CV$-category $[\CC^{\mathrm{op}},\CV]$ and its $\CV$-category 
$[\CC^{\mathrm{op}},\CV]^{\widebar{\Gamma}}$ of
Eilenberg-Moore $\widebar{\Gamma}$-coalgebras.

The following definition agrees with the category $\CC^{\Gamma}$ defined by Thi\'ebaud \cite{Thi\'ebaudPhD} in the case $\CV = \mathrm{Set}$.
   
\begin{definition}\label{Gammaalgebras}
{\em The $\CV$-category of $\Gamma$-algebras in $\CC$} is defined by the pullback \eqref{Gamma-alg} in $\CV\text{-}\mathrm{Cat}$ of the underlying $\CV$-functor along the Yoneda embedding $\yo$.
  \begin{eqnarray}\label{Gamma-alg}
  \begin{aligned}
\xymatrix{
\CC^{\Gamma} \ar[rr]^-{\yo^{\Gamma}} \ar[d]_-{\mathrm{und}} & & [\CC^{\mathrm{op}},\CV]^{\widebar{\Gamma}} \ar[d]^-{\mathrm{und}} \\
\CC \ar[rr]_-{\yo} & & [\CC^{\mathrm{op}},\CV]}
\end{aligned}
\end{eqnarray}
So such a $\Gamma$-algebra consists of an object $C \in \CC$ equipped with a {\em coaction} morphism $\gamma : I \to \Gamma(C,C)$, subject to the two axioms \eqref{axx: Gamma-alg}.
We will write $\gamma_X : \CC(C,X) \to \Gamma(C,X)$ for the natural family corresponding to
$\gamma$ under the Yoneda bijection. Similarly we have $\gamma^X : \CC(X,C) \to \Gamma(X,C)$.  
 \begin{eqnarray}\label{axx: Gamma-alg}
  \begin{aligned}
\xymatrix{
I \ar[d]_-{\gamma} \ar[rd]^-{1_C} & \\
\Gamma(C,C) \ar[r]_-{\varepsilon_{C,C}} & \CC(C,C)}
\qquad
\xymatrix{
I \ar[r]^-{\gamma\ox\gamma} \ar[d]_-{\gamma} & \Gamma(C,C)\ox \Gamma(C,C)\ar[d]^-{\mathrm{in}_C} \\
\Gamma(C,C) \ar[r]_-{\delta_{C,C}} & \int^X \Gamma(X,C) \ox \Gamma(C,X)}
\end{aligned}
\end{eqnarray}
We will call any $\CV$-functor into $\CC$, isomorphic over $\CC$ to $\mathrm{und} : \CC^{\Gamma}\to \CC$, {\em Thi\'ebaud algebraic over $\CC$}.  
\end{definition}

\begin{example}\label{exxprocomonad} The construction $\CC^{\Gamma}$ includes the Eilenberg-Moore
constructions for both $\CV$-monads and $\CV$-comonads.
\begin{itemize}
\item[1.] If $T = (T,\eta,\mu)$ is a $\CV$-monad on the $\CV$-category $\CC$ and we take $\Gamma = T^*$ so that $\Gamma(X,Y)=\CC(TX,Y)$ with counit $\varepsilon$ and comultiplication $\delta$ induced by the unit $\eta$ and multiplication $\mu$
then $\CC^{\Gamma} \cong \CC^T$, the $\CV$-category of $T$-algebras.  
  
\item[2.] If $G = (G,\varepsilon,\delta)$ is a $\CV$-comonad on the $\CV$-category $\CC$ and we take $\Gamma = G_*$ so that $\Gamma(Y,Z)=\CC(Y,GZ)$ with
counit and comultiplication induced by those of $G$ then $\CC^{\Gamma} \cong \CC^G$, the $\CV$-category of $G$-coalgebras.   
\end{itemize}
\end{example}

The two main closure properties Thi\'ebaud proved in \cite{Thi\'ebaudPhD} were that Thi\'ebaud algebraicity is closed under pullback and exponentiation. We now look at that. 

Given a $\CV$-functor $W : \CD\to \CC$ and a $\CV$-procomonad $\Gamma$ on $\CC$, we have the $\CV$-procomonad $\Gamma_W = W^*\circ \Gamma \circ W_*$
on $\CD$; it is {\em the lifting of $\Gamma$ through $W_*$} in $\mathfrak{M}$ (see Section 2 of \cite{3}). 

\begin{proposition}\label{pbGammaAlgs} The following square is a pullback.
 \begin{eqnarray*}
  \begin{aligned}
\xymatrix{
\CD^{\Gamma_W} \ar[rr]^-{ } \ar[d]_-{\mathrm{und}} & & \CC^{\Gamma} \ar[d]^-{\mathrm{und}} \\
\CD \ar[rr]_-{W} & & \CC}
\end{aligned}
\end{eqnarray*}
\end{proposition}
\begin{proof}
Using Yoneda, we deduce that $\Gamma_W(D',D) \cong \Gamma(WD',WD)$.
The remaining details are routine.    
\end{proof}
\begin{corollary} Thi\'ebaud algebraicity is the closure under pullback of comonadicity. 
\end{corollary}

Given $\CV$-categories $\CA$ and $\CC$, each $\CV$-procomonad $\Gamma$ on $\CC$
defines a $\CV$-procomonad $\Gamma^{\CA}$ on $[\CA,\CC]$ via the commutative diagram \eqref{Gammasuper}; it is the lifting of $[1_{\CA},\widebar{\Gamma}]_*$ through
$[1_{\CA},\yo]_*$ and satisfies the simple formula:
\begin{eqnarray*}
\Gamma^{\CA}(F',F) = \int_A\Gamma(F'A,FA) \ .
\end{eqnarray*}
 
\begin{eqnarray}\label{Gammasuper}
\begin{aligned}
\xymatrix{
[\CA,\CC] \ar[rr]^-{\Gamma^{\CA}} \ar[d]_-{[1_{CA},\yo]_*} && [\CA,\CC]  \\
[\CA,[\CC^{\mathrm{op}},\CV]] \ar[rr]_-{[1_{\CA},\widebar{\Gamma}]_*} && [\CA,[\CC^{\mathrm{op}},\CV]]
\ar[u]_-{[1_{\CA},\yo]^*}}
\end{aligned}
\end{eqnarray}

\begin{proposition}\label{expGammaAlgs} $[\CA,\CC]^{\Gamma^{\CA}} \cong [\CA,\CC^{\Gamma}]$ over $[\CA,\CC]$.
\end{proposition}
\begin{proof}
The pullback \eqref{Gamma-alg} is preserved by exponentiation $[\CA,-]$ by $\CA$. 
The Eilenberg-Moore construction of coalgebras in $\CV\text{-}\mathrm{Cat}$ is also
preserved by exponentiation. So we have the pullback
\begin{eqnarray*}
\xymatrix{
[\CA,\CC^{\Gamma}] \ar[rr]^-{ } \ar[d]_-{[1_{\CA},\mathrm{und}]} & & [\CA, [\CC^{\mathrm{op}},\CV]]^{[1_{\CA},\widebar{\Gamma}]} \ar[d]^-{\mathrm{und}} \\
[\CA,\CC] \ar[rr]_-{[1_{\CA},\yo]} & & [\CA,[\CC^{\mathrm{op}},\CV]] \ .}
\end{eqnarray*}
The result now follows from Proposition~\ref{pbGammaAlgs}, the definition \eqref{Gammasuper} of $\Gamma^{\CA}$, and the second of Example~\ref{exxprocomonad}.   
\end{proof}

An object $\CX$ of a monoidal bicategory $\mathfrak{N}$ (see \cite{60}) is {\em magmal} when a 1-morphism $P : \CX\ox \CX\to \CX$ is specified. For example, every monoidale (called pseudomonoid by \cite{60} and elsewhere) in $\mathfrak{N}$ has an underlying magmal object. 

In particular, a magmal object in $\CV\text{-}\mathrm{Cat}$ is called a magmal $\CV$-category
as in Definition~\ref{defmag}. 
A module $M : \CC\nrightarrow \CD$ between magmal $\CV$-categories is called 
{\em magmal} when it is equipped with a $\CV$-natural family $M_2$ of morphisms 
\begin{eqnarray*}
M_{2; C,C'}^{\ D,D'} : M(D,C)\ox M(D',C') \to M(D\ox D',C\ox C') \ ;
\end{eqnarray*}
such families, by the universal property of coend and Yoneda's Lemma, are in bijection with
$\CV$-natural families of morphisms 
\begin{eqnarray*}
M_{2; C,C'}^{\ D''} : \int^{D,D'} M(D,C)\ox M(D',C')\ox \CD(D'',D\ox D') \to M(D'',C\ox C') \ ;
\end{eqnarray*}
and in bijection with $\CV$-natural families of morphisms 
\begin{eqnarray*}
M_{2; C''}^{\ D,D'} : \int^{C,C'}\CC(C\ox C',C'')\ox M(D,C)\ox M(D',C') \to M(D\ox D',C'') \ .
\end{eqnarray*}

A module morphism $\alpha : M\Ra N$ is {\em magmal} when 
$$\alpha_{D\ox D',C\ox C'}\circ M_{2; C,C'}^{\ D,D'} = N_{2; C,C'}^{\ D,D'} \circ (\alpha_{D,C}\ox \alpha_{D',C'}) \ .$$ 
 
A $\CV$-functor $S : \CC\to \CD$ is magmal as in Definition~\ref{defmag} if and only if the module $S_*$ is. 
Using the Yoneda Lemma, we see that this amounts to a $\CV$-natural family of morphisms    
\begin{eqnarray*}
S_2 = S_{2; C,C'} : SC\ox SC' \lra S(C\ox C') 
\end{eqnarray*} 
as in diagram \eqref{magfun}.
Call $S$ {\em strong magmal} when all $S_{2, C,C'}$ are invertible. 
 
If $\CC$ is a small magmal $\CV$-category then the presheaf $\CV$-category $\widehat{\CC} = [\CC^{\mathrm{op}},\CV]$ has the
{\em Day convolution magmal structure} 
$$\asterisk : \widehat{\CC} \ox \widehat{\CC} \to \widehat{\CC}$$ 
defined by 
\begin{eqnarray*}
(F\asterisk F')Z = \int^{X,Y}\CC(Z,X\ox Y)\ox FX\ox F'Y \ .
\end{eqnarray*}
The Yoneda embedding $\yo : \CC\to \widehat{\CC}$ is strong magmal and
has the (bicategorical) universal property of the magmal small cocompletion of $\CC$:
for any small-cocomplete magmal $\CX$, the category of 
[strong-]magmal $\CV$-functors $\CC\to \CX$ is equivalent 
(via left Kan extension along $\yo$) to the category of colimit-preserving
[strong-]magmal $\CV$-functors $\widehat{\CC} \to \CX$. 

A procomonad $\Gamma = (\Gamma,\varepsilon , \delta)$ on a $\CV$-category $\CX$ is {\em magmal} when $\Gamma,\varepsilon , \delta$ are all magmal.

Let $\Gamma$ be a magmal procomonad on the magmal $\CV$-category $\CC$.
Then $\widebar{\Gamma}$ is a magmal comonad on $\widehat{\CC}$.
To obtain $\widebar{\Gamma}_2$ we use the isomorphisms
\begin{eqnarray*}
\widebar{\Gamma}(F)\asterisk \widebar{\Gamma}(F') \cong \int^{UVXY}FU\ox F'V\ox \Gamma(X,U)\ox \Gamma(Y,V) \ox \CC(-,X\ox Y) 
\end{eqnarray*}
and
\begin{eqnarray*}
\widebar{\Gamma}(F\asterisk F') \cong \int^{UV}FU\ox F'V\ox \Gamma(-,U\ox V) \ , 
\end{eqnarray*}
to transport $\int^{UV}1_{FU}\ox 1_{F'V}\ox \Gamma_{2; U,V}^{\ -}$ to obtain $\widebar{\Gamma}_{2; F, F'} : \widebar{\Gamma}(F)\asterisk \widebar{\Gamma}(F')\to \widebar{\Gamma}(F\asterisk F')$.
Using \eqref{magcoalg}, we obtain a magmal structure on $\widehat{\CC}^{\ \widebar{\Gamma}}$. This restricts along the fully faithful $\CV$-functor
$\yo^{\Gamma} : \CC^{\Gamma} \to \widehat{\CC}^{\ \widebar{\Gamma}}$ 
to a magmal structure on $\CC^{\Gamma}$ which is defined by
\begin{eqnarray*}
\lefteqn{\left(X, \ I\xra{\xi}\Gamma(X,X)\right)\ox \left(Y, \ I\xra{\upsilon}\Gamma(Y,Y)\right)} \\
& : = & \left(X\ox Y, \ I\xra{\xi\ox \upsilon}\Gamma(X,X)\ox \Gamma(Y,Y)\xra{\Gamma_{2;XY}^{XY}}\Gamma(X\ox Y,X\ox Y)\right) \ . 
\end{eqnarray*}

%
\begin{definition}\label{strongcloakal}
Suppose $S : \CC\to \CD$ is a $\CV$-functor between magmal $\CV$-categories.
Suppose $Z$ is left cloaked by $Y$ in $\CC$. 
We say {\em $S$ preserves the left cloaking of $Z$ by $Y$} when $S[Y,Z]$ provides
a left cloaking of $SZ$ by $SY$. If this holds for all $Y, Z$ then we say $S$ is 
{\em strong left cloakal}; by Proposition~\ref{magmalcloakal}, it follows that $S$
is magmal, but it is not necessarily strong magmal.        
\end{definition}

\begin{proposition}\label{cloakalpresheaves}
For any magmal $\CV$-category $\CC$, the convolution magmal presheaf $\CV$-category $\widehat{\CC}$ is left cloakal. For $H, K\in \widehat{\CC}$, the left cloaking of $K$ by $H$ is given by 
$$[H,K]U = \int_V[HV,K(U\ox V)] \cong \ \widehat{\CC}(H,K(U\ox -)) \ .$$
In particular, $[\yo Y, \yo Z]X\cong \CC(X\ox Y,Z)$ so that the Yoneda embedding 
$\yo : \CC\to \widehat{\CC}$ preserves any left cloakings $\CC$ admits. 
\end{proposition}
\begin{proof}
This follows {\em mutatis mutandis} the proof by Day \cite{DayConv} in the monoidal case.  
\end{proof}

\begin{definition}\label{fusionGamma} The {\em fusion morphisms for magmal procomonad $\Gamma$} are 
the Wood fusion morphisms for the magmal comonad $\widebar{\Gamma}$ restricted to representables.
\end{definition}

Let us make this definition more explicit. According to Definition~\ref{defwoodfusion}, Wood fusion for $\widebar{\Gamma}$ is the composite  
\begin{eqnarray}\label{watrhoK}
\mathrm{w}_{\rho, K} : = \widebar{\Gamma}[H,K] \xra{\widebar{\Gamma}^{\ell}_2} [\widebar{\Gamma}H,\widebar{\Gamma}K] \xra{[\rho,1]} [H,\widebar{\Gamma}K] \ ,
\end{eqnarray} 
for $(H,\rho)\in \widehat{\CC}^{\ \widebar{\Gamma}}$ and $K\in \widehat{\CC}$. 
For $X, Z\in \CC$ and $(Y,\upsilon)\in \CC^{\Gamma}$, put $K = \yo Z$, $H = \yo Y$, and $\rho = \upsilon^-$,
then we write $\mathrm{w}_{X, \upsilon, Z}$ instead of $(\mathrm{w}_{\rho, K})_X$.   
Using Yoneda, we obtain $$(\widebar{\Gamma}\yo Y)X \cong \Gamma(X,Y) \ ,$$
$$(\widebar{\Gamma}[\yo Y,\yo Z])X \cong \int^U\CC(U\ox Y,Z)\ox \Gamma(X,U) \ ,$$ 
$$[\widebar{\Gamma}\yo Y,\widebar{\Gamma}\yo Z]X = \widehat{\CC}(\Gamma(-,Y),\Gamma(X\ox -,Z)) \ , $$ 
$$[\yo Y,\widebar{\Gamma}\yo Z]X = \Gamma(X\ox Y,Z) \ . $$

\begin{proposition} Let $\Gamma$ be a magmal procomonad on the magmal $\CV$-category $\CC$. Suppose $(Y,I \xra{\upsilon}\Gamma(Y,Y))\in \CC^{\Gamma}$.
The diagram \eqref{magprocomonadfusion} commutes.
{\footnotesize
\begin{eqnarray}\label{magprocomonadfusion}
\begin{aligned}
\xymatrix{
\int^{U}\CC(U\ox Y,Z)\ox\Gamma(X,U) \ar[r]^-{\mathrm{w}_{X\upsilon Z}} \ar[d]_-{\int^{U} 11\upsilon}  & \Gamma(X\ox Y,Z) \\
  \int^{U}\CC(U\ox Y,Z)\ox \Gamma(X,U)\ox \Gamma(Y,Y) \ar[r]_-{\int^{U}\mathrm{in}_{Y}}   &   \int^{UV}\CC(U\ox V,Z)  \ox\Gamma(X,U)\ox\Gamma(Y,V) \ar[u]_-{\Gamma_{2;Z}^{XY}} 
  }
\end{aligned} 
\end{eqnarray} 
}
\end{proposition}
\begin{proof} Replace the end vertex $\Gamma(X\ox Y,Z)$ of the diagram by the Yoneda isomorphic $\widehat{\CC}(\CC(-,Y),\Gamma(X\ox -, Z))$. The morphism $\mathrm{w}_{X\upsilon Z}$ transports to the composite {\footnotesize
$$ \int^U\CC(U\ox Y,Z)\ox \Gamma(X,U) \xra{\widebar{\Gamma}^{\ell}_2} \widehat{\CC}(\Gamma(-,Y),\Gamma(X\ox -,Z))\xra{\widehat{\CC}(\upsilon^-, 1)} \widehat{\CC}(\CC(-Y),\Gamma(X\ox -, Z)) \ .$$} 
It suffices to show that the two paths around \eqref{magprocomonadfusion} agree after we precompose the diagram with each injection 
$$\CC(U\ox Y,Z)\ox \Gamma(X,U) \xra{\mathrm{in}_U} \int^{U}\CC(U\ox Y,Z)\ox \Gamma(X,U)$$
and postcompose with each projection
$$\widehat{\CC}(\CC(-,Y),\Gamma(X\ox -, Z))\xra{\mathrm{pr}_V} [\CC(V,Y),\Gamma(X\ox V, Z)] \ .$$
Now we need to show that two paths $$\CC(U\ox Y,Z)\ox \Gamma(X,U) \to [\CC(V,Y),\Gamma(X\ox V, Z)]$$
are equal. By Yoneda, it suffices to check them equal after taking $Z = U\ox Y$ and evaluating at the identity
(that is, on precomposing with 
$$j_{U\ox Y}\ox 1_{\Gamma(X,U)}: \Gamma(X,U)\to \CC(U\ox Y,U\ox Y)\ox \Gamma(X,U) \ ) \ .$$
Two commutative diagrams then show that both paths reduce to the morphism
$$\Gamma(X,U) \lra [\CC(V,Y),\Gamma(X\ox V, U\ox Y)]$$
corresponding to
$$\Gamma(X,U) \ox \CC(V,Y)\xra{1\ox \upsilon^V} \Gamma(X,U)\ox \Gamma(V,Y)\xra{\Gamma_2} \Gamma(X\ox V, U\ox Y)$$
under the closed-monoidal adjunction $-\ox\CC(V,Y)\dashv [\CC(V,Y), -]$ for $\CV$.   
\end{proof}

\begin{corollary}\label{cloakedfusion}
Suppose $Z$ is left cloaked by $Y$ in the magmal $\CV$-category $\CC$ and $(Y,\upsilon)\in \CC^{\Gamma}$.
Then the fusion morphism \eqref{magprocomonadfusion} becomes the composite \eqref{cloakedfusion}.
\begin{eqnarray}\label{cloakedfusion}
\begin{aligned}
\xymatrix{
\Gamma(X,[Y,Z]) \ar[r]^-{\mathrm{w}_{X\upsilon Z}} \ar[d]_-{1\ox\upsilon}  & \Gamma(X\ox Y,Z) \\
  \Gamma(X,[Y,Z]) \ox \Gamma(Y,Y) \ar[r]_-{\Gamma_2}   &  \Gamma(X\ox Y,[Y,Z]\ox Y) \ar[u]_-{\Gamma(1,\mathrm{ev})} 
  }
\end{aligned} 
\end{eqnarray}
\end{corollary}
\begin{proof}
We leave this as an exercise for the reader.
\end{proof}

\begin{definition} We call $\Gamma$ {\em Hopf at $(Y,\upsilon)$} when $\mathrm{w}_{X\upsilon Z}$ is invertible for all $X, Z$.
\end{definition}

\begin{lemma}\label{FfromZ} 
Suppose $(Y,\upsilon)\in \CC^{\Gamma}$ is such that $[Y,Z]\in \CC$ exists for all $Z\in \CC$.
Then $\Gamma$ is Hopf at $(Y,\upsilon)$ if and only if $\widebar{\Gamma}$ is Hopf at $\yo^{\Gamma}(Y,\upsilon)$.  
\end{lemma}
\begin{proof} ``If'' is clear since $\mathrm{w}_{X\upsilon Z}$ is defined as a special case of Wood fusion for $\widebar{\Gamma}$.
Conversely, suppose $\mathrm{w}_{X\upsilon Z} : \Gamma(X,[Y,Z])\cong \Gamma(X\ox Y,Z)$ for all $X,Z\in \CC$.
For $K\in \widehat{\CC}$, one easily calculates that $[\yo Y,K]\cong K(-\ox Y)$, so
\begin{eqnarray*}
\widebar{\Gamma}[\yo Y,K] & \cong &  \int^UK(U\ox Y)\ox \Gamma(-,U) \\
& \cong & \int^{UZ}\CC(U\ox Y,Z)\ox KZ\ox \Gamma(-,U) \\
& \cong & \int^{UZ}\CC(U, [Y,Z])\ox KZ\ox \Gamma(-,U) \\
& \cong & \int^{Z} KZ\ox \Gamma(-,[Y,Z]) \\
& \cong &  \int^{Z} KZ\ox \Gamma(-\ox Y,Z) \\
& \cong & \widebar{\Gamma}(K)(-\ox Y) \\
& \cong & [\yo Y, \widebar{\Gamma}(K)] \ ,
\end{eqnarray*}  
where the fifth isomorphism uses the invertible $\Gamma$-fusion morphisms $\mathrm{w}_{-\upsilon Z}$. 
The composite is \eqref{watrhoK} for $H = \yo Y$ and $\rho = \upsilon^-$.
\end{proof}

Suppose $Z$ is left cloaked by $Y$ in the magmal $\CV$-category $\CC$ and $(Y,\upsilon), (Z,\zeta) \in \CC^{\Gamma}$. Suppose $\Gamma$ is Hopf at $(Y,\upsilon)$. Then we have $([Y,Z],\omega)\in \CC^{\Gamma}$ defined by commutativity of \eqref{omega}. 
\begin{eqnarray}\label{omega}
\begin{aligned}
\xymatrix{
I \ar[r]^-{\omega} \ar[d]_-{\zeta}  & \Gamma([Y,Z],[Y,Z]) \ar[d]^-{\mathrm{w}_{[Y,Z]\upsilon Z}} \\
  \Gamma(Z,Z)  \ar[r]_-{\Gamma(\mathrm{ev},1)}   &  \Gamma([Y,Z]\ox Y,Z)  
  }
\end{aligned} 
\end{eqnarray}

Via Example~\ref{exxprocomonad}, we obtain our unification of \cite{Bruguieres2011} and \cite{Wood1978}
as an application of Proposition~\ref{magcomoncloaks} with $G = \widebar{\Gamma}$, using Lemmas~\ref{pbcloakcreation} and \ref{FfromZ}. 

\begin{theorem}\label{magprocomoncloaks} Let $ {\Gamma}$ be a magmal procomonad on a magmal $\CV$-category $\CC$.
Suppose $(Y,\upsilon)\in \CC^{ \Gamma}$ is such that $[Y,Z]\in \CC$ exists for all $Z\in \CC$. 
Then $\Gamma$ is Hopf at $(Y,\upsilon)$ if and only if $\mathrm{und} : \CC^{ \Gamma} \to \CC$ creates all cloaks by $(Y,\upsilon)$.
In this case, for any $(Z,\zeta)\in \CC^{ {G}}$, 
$$[(Y,\upsilon),(Z,\zeta)] \cong ([Y,Z], \omega) \ ,$$
where $\omega$ is defined by \eqref{omega}. 
\end{theorem}

\begin{remark} In general the $\CV$-functor $\mathrm{und} : \CC^{ \Gamma} \to \CC$ has neither left nor right adjoint.
Example~\ref{exxprocomonad} is where it does.
 
\end{remark}
 \appendix
 
 \section{Lifting adjunctions and doctrinal adjunction}\label{Lada}

After my talk on Wood fusion in the Australian Category Seminar on 8 February 2023, Steve Lack and
Richard Garner suggested that the results I presented were obtainable from adjoint lifting theorems and that I should 
look at Peter Johnstone's paper \cite{Johnstone1975}. This section addresses that suggestion.

Richard Wood \cite{Wood1978} already referred
to William Keigher \cite{Keigher1975}. Johnstone states he learned of \cite{Keigher1975} after writing \cite{Johnstone1975}. 

We revisit the Adjoint Triangle Theorem of Eduardo Dubuc \cite{DubucAT} in Appendix~\ref{liftings}. 
Now we will see that other adjoint lifting results can be viewed as consequences of doctrinal
adjunction (in the sense of Max Kelly in \cite{Kelly1974}) involving examples as in Theorem 9 of \cite{3} 
and Theorem 1 of \cite{4}.

Let $\mathrm{Mnd}\mathfrak{C}$ denote the bicategory of monads in a bicategory $\mathfrak{C}$ essentially as defined in \cite{3}.
An object is a pair $(A,s)$ consisting of an object $A\in \mathfrak{C}$ and a monoid $s$ (called a {\em monad on $A$}) 
in the endomonoidal category $\mathfrak{C}(A,A)$. The unit and multiplication of $s$ will be denoted by
$\eta : 1_A\Rightarrow s$ and $\mu : s\circ s\Rightarrow s$. A morphism $(u,\phi) : (A,s)\to (B,t)$ 
(called a {\em monad morphism}) consists of a morphism $u : A\to B$ equipped with a 2-morphism
\begin{equation}\label{rif}
 \begin{aligned}
  \xymatrix{
A \ar[d]_{u}^(0.5){\phantom{aaa}}="1" \ar[rr]^{s}  && A \ar[d]^{u}_(0.5){\phantom{aaa}}="2" \ar@{=>}"1";"2"^-{\phi}
\\
B \ar[rr]_-{t} && B 
}
 \end{aligned}
\end{equation}
in $\mathfrak{C}$ compatible with $\eta$, $\mu$ in the obvious way.
A 2-morphism $\sigma : (u,\phi)\Rightarrow (v,\psi) : (A,s)\to (B,t)$ is a 2-morphism $\sigma : u \Rightarrow v$ in $\mathfrak{C}$ such that $\sigma s\circ \phi = \psi \circ t\sigma$.
Composition is performed by pasting. There is a forgetful pseudofunctor 
\begin{equation}\label{MndKtoK}
\mathrm{Mnd}\mathfrak{C}\to \mathfrak{C} \ , \ (A,s) \mapsto A \ .
\end{equation}
A morphism in $\mathrm{Mnd}^{\mathrm{op}}\mathfrak{C} = (\mathrm{Mnd}\mathfrak{C}^{\mathrm{op}})^{\mathrm{op}}$ 
is called a {\em monad opmorphism}: the 2-morphism in \eqref{rif} is reversed.
A morphism in $\mathrm{Mnd}^{\mathrm{co}}\mathfrak{C} = (\mathrm{Mnd}\mathfrak{C}^{\mathrm{co}})^{\mathrm{co}}$ 
is called a {\em comonad morphism}. A morphism in $\mathrm{Mnd}^{\mathrm{coop}}\mathfrak{C} = (\mathrm{Mnd}\mathfrak{C}^{\mathrm{coop}})^{\mathrm{coop}}$ 
is called a {\em comonad opmorphism}.
\begin{example}\label{fusionmndmor}
\begin{itemize} 
\item[1.] Let $T$ be an opmagmal monad on the magmal $\CV$-category $\CC$ as in Section~\ref{Ffom}. For $(X,\alpha)\in \CC^T$, $T$-fusion supplies a monad morphism
\begin{equation*}
 \begin{aligned}
  \xymatrix{
\CC \ar[d]_{-\ox X}^(0.5){\phantom{aaa}}="1" \ar[rr]^{T}  && \CC \ar[d]^{-\ox X}_(0.5){\phantom{aaa}}="2" \ar@{=>}"1";"2"^-{v_{-,\alpha}}
\\
\CC \ar[rr]_-{T} && \CC \ . 
}
 \end{aligned}
\end{equation*}
\item[2.] Let $G$ be a magmal comonad on the magmal $\CV$-category $\CC$ as in Definition~\ref{defmag}. For $(Y,\upsilon)\in \CC^G$, Wood fusion supplies a comonad opmorphism
\begin{equation*}
 \begin{aligned}
  \xymatrix{
\CC \ar[d]_{[Y,-]}^(0.5){\phantom{aaa}}="1" \ar[rr]^{G}  && \CC \ar[d]^{[Y,-]}_(0.5){\phantom{aaa}}="2" \ar@{=>}"1";"2"^-{w_{-,\upsilon}}
\\
\CC \ar[rr]_-{G} && \CC \ . 
}
 \end{aligned}
\end{equation*}
\end{itemize}
\end{example}
Part of doctrinal adjunction is the fact that, if $(u,\phi) : (A,s)\to (B,t)$ is a monad morphism and $f \dashv u$
is an adjunction, then $(f,\hat{\phi}) : (B,t)\to (A,s)$ is a monad opmorphism where $\hat{\phi} : ft\Rightarrow sf$ is the mate of $\phi : tu\Rightarrow us$. The other part is obtained by examining adjunctions      
\begin{equation}\label{mndtocat}
(f,\theta)\dashv (u,\phi) : (A,s)\to (B,t)
\end{equation}
in $\mathrm{Mnd}\mathfrak{C}$. Since pseudofunctors preserve adjunctions, we use \eqref{MndKtoK} to deduce that
$f\dashv u$ in $\mathfrak{C}$ and that the counit $\alpha : fu\Rightarrow 1_A$ and unit $\beta : 1\Rightarrow uf$ must 
be 2-morphisms in $\mathrm{Mnd}\mathfrak{C}$. Using only the first of these, a little diagram shows that $\theta$
has inverse the mate 
\begin{equation}\label{mateofphi}
\hat{\phi} : ft\xra{ft\beta} ftuf \xra{f\phi f} fusf \xra{\alpha sf} sf
\end{equation} 
of $\phi$. Yet, if this mate of $\phi$ has an inverse at all, one sees that both $\alpha$ and $\beta$ are 2-morphisms in $\mathrm{Mnd}\mathfrak{C}$. This proves:

\begin{proposition}\label{doctadjMnd}
A morphism $(u,\phi) : (A,s)\to (B,t)$ has a left adjoint in $\mathrm{Mnd}\mathfrak{C}$ if and only if $u : A\to B$
has a left adjoint and the mate \eqref{mateofphi} of $\phi$ is invertible.
\end{proposition}
 
Let $\mathrm{Fun}\mathfrak{C}$ denote the lax morphism bicategory of $\mathfrak{C}$. 
An object is a morphism $X\xra{x} A$ in $\mathfrak{C}$. A morphism $(u,\upsilon, \bar{u}) : x\to y$ is a diagram 
\begin{equation}\label{morfun}
 \begin{aligned}
  \xymatrix{
X \ar[d]_{\bar{u}}^(0.5){\phantom{aaa}}="1" \ar[rr]^{x}  && A \ar[d]^{u}_(0.5){\phantom{aaa}}="2" \ar@{=>}"1";"2"^-{\upsilon}
\\
Y \ar[rr]_-{y} && B 
}
 \end{aligned}
\end{equation}
in $\mathfrak{C}$.
A 2-morphism $(\sigma, \bar{\sigma}) : (u,\upsilon, \bar{u})\Rightarrow (v,\omega, \bar{v}) : x\to y$
is a pair of 2-morphisms in $\mathfrak{C}$ satisfying $\sigma x \circ \upsilon = \omega \circ y\bar{\sigma}$. 
A morphism $(u,\upsilon, \bar{u})$ is called {\em strong} when $\upsilon$ is invertible. Let $\mathrm{sFun}\mathfrak{C}$ denote the sub-2-category of $\mathrm{Fun}\mathfrak{C}$ obtained by restricting to the strong morphisms.   

Let us look at adjunctions $(f,\tau, \bar{f}) \dashv (u,\upsilon, \bar{u})$ in $\mathrm{Fun}\mathfrak{C}$.
As before, because of the existence of forgetful pseudofunctors, we must have adjunctions $f\dashv u$ and $\bar{f}\dashv \bar{u}$
in $\mathfrak{C}$ such that the counits and units form 2-morphisms $(\alpha, \bar{\alpha})$ and $(\beta, \bar{\beta})$
in $\mathrm{Fun}\mathfrak{C}$. The first of these yields that the mate $\hat{\upsilon}$ of $\upsilon$ is a left inverse
for $\tau$ while the second yields that $\hat{\upsilon}$ is a right inverse for $\tau$. On the other hand,
any inverse for $\hat{\upsilon}$ does render $(\alpha, \bar{\alpha})$ and $(\beta, \bar{\beta})$ 2-morphisms.      

\begin{proposition}\label{doctadjFun}
A morphism $(u,\upsilon, \bar{u}) : x\to y$ has a left adjoint in $\mathrm{Fun}\mathfrak{C}$ if and only if both 
$u : A\to B$ and $\bar{u} : X\to Y$ have left adjoints and the mate $\hat{\upsilon}$ of $\upsilon$ is invertible.
Any left adjoint in $\mathrm{Fun}\mathfrak{C}$ is in $\mathrm{sFun}\mathfrak{C}$. Any $(f,\tau, \bar{f})$ in $\mathrm{sFun}\mathfrak{C}$ has a right adjoint in $\mathrm{Fun}\mathfrak{C}$ if and only if $f$ and $\bar{f}$ have right adjoints in $\mathfrak{C}$.   
\end{proposition}

Now suppose $\mathfrak{C}$ admits the construction of algebras in the bicategorical sense: for each monad $(A,s)$,
there is an Eilenberg-Moore $s$-algebra $x_s : A^s \to A$ with action $\xi_s : sx_s \to x_s$ for which the functor
\begin{eqnarray*}
\mathfrak{C}(X,A^s) \lra \mathfrak{C}(X,A)^{\mathfrak{C}(1_X,s)} \ , \ (X\xra{h} A^s) \mapsto (x_sh, s x_s h\xra{\xi_s h}x_sh)
\end{eqnarray*}
is an equivalence for all $X\in \mathfrak{C}$. 
Then we have a pseudofunctor
\begin{equation}\label{em}
\mathrm{EM} : \mathrm{Mnd}\mathfrak{C}\to \mathrm{sFun}\mathfrak{C}
\end{equation}
defined as follows. For each monad $(A,s)$, we put $\mathrm{EM}(A,s) = x_s$. 
For each monad morphism $(u,\phi) : (A,s)\to (B,t)$, we have a $t$-algebra
\begin{eqnarray}\label{talgfromphi}
tux_s \xra{\upsilon x_s} usx_s \xra{u\xi_s}ux_s \ ,
\end{eqnarray}
so there exist (uniquely up to isomorphism) a morphism $\bar{u} : A^s\to A^t$ and an isomorphism $x_t\bar{u}\xra{\upsilon} ux_s$ such that $\upsilon$ becomes a $t$-algebra isomorphism from
$(x_t\bar{u},\xi_t\bar{u})$ to the $t$-algebra \eqref{talgfromphi}. 
We put $\mathrm{EM}(u,\phi) = (u, \upsilon, \bar{u}) : x_s\to x_t$; it is a strong morphism.
For a 2-morphism $\sigma : (u,\phi)\Rightarrow (v,\psi) : (A,s)\to (B,t)$, there is a unique
2-morphism $\bar{\sigma} : \bar{u}\Rightarrow \bar{v}$ such that
$(\sigma, \bar{\sigma}) : (u,\upsilon, \bar{u})\Rightarrow (v,\omega, \bar{v}) : x_s\to x_t$
is a 2-morphism in $\mathrm{Fun}\mathfrak{C}$; we put $\mathrm{EM}\sigma = (\sigma, \bar{\sigma})$.

Here is a restatement of an observation of Appelgate \cite{Appel1965}; also see Lemma 1 of \cite{Johnstone1975}.
\begin{proposition}\label{Appel}
The pseudofunctor \eqref{em} is an equivalence on homcategories. 
\end{proposition}
\begin{proof}
Take $(u, \upsilon, \bar{u}) : x_s\to x_t$ in $\mathrm{sFun}\mathfrak{C}$ and let $\tau : y_tu\to \bar{u}y_s$
be the mate of $\upsilon^{-1} : ux_s\to x_t\bar{u}$ under the adjunctions $y_s\dashv x_s$ and
$y_t\dashv x_t$. Now put
\begin{eqnarray}\label{locequivEM}
\phi  = (tu\cong x_ty_tu\xra{x_t\tau} x_t\bar{u}y_s\xra{\upsilon y_s} ux_sy_s \cong us) \ . 
\end{eqnarray}
Recalling how the unit and multiplication of the monads are obtained from the unit and counit of the generating adjunctions, we routinely check that $(u,\phi) : (A,s)\to (B,t)$ is a morphism of $\mathrm{Mnd}\mathfrak{C}$ (a string diagram proof is attractive) and that $\mathrm{EM}(u,\phi) \cong (u, \upsilon, \bar{u})$.
From the definition of $\mathrm{EM}$ on 2-morphisms we see, for each $(\sigma, \bar{\sigma}) : \mathrm{EM}(u,\phi)\Rightarrow \mathrm{EM}(v,\psi) : x_s\to x_t$, that $\sigma : (u,\phi)\Rightarrow (v,\psi)$ is the unique 2-morphism of $\mathrm{Mnd}\mathfrak{C}$ with
$\mathrm{EM}\sigma = (\sigma, \bar{\sigma})$.            
\end{proof}

\begin{remark}
In the above proof, notice that, since composing with $x_t$ is conservative, $\phi$ is invertible if and only if $\tau$ is. This is Lemma 3 of Johnstone \cite{Johnstone1975}.
\end{remark}

\begin{corollary}\label{EMadjprop}
A monad morphism $(u,\phi) : (A,s)\to (B,t)$ has a left adjoint in $\mathrm{Mnd}\mathfrak{C}$ if and only if
$\mathrm{EM}(u,\phi)$ has a left adjoint in $\mathrm{Fun}\mathfrak{C}$. 
\end{corollary}

The second sentence of Corollary~\ref{PTJKeig} is Theorem 4 of \cite{Johnstone1975} and the dual
of Corollary 2.3 of \cite{Keigher1975}. Furthermore, by way of Example~\ref{fusionmndmor}, 
it relates to Theorem 3.6 of \cite{Bruguieres2011} and to our Proposition~\ref{magcomoncloaks}.
 
\begin{corollary}\label{PTJKeig}
A monad morphism $(u,\phi) : (A,s)\to (B,t)$ has a right adjoint in $\mathrm{Mnd}\mathfrak{C}$ if and only if
$\mathrm{EM}(u,\phi)$ has a right adjoint $(r,\rho, \bar{r})$ in $\mathrm{sFun}\mathfrak{C}$. In particular,
if $u$ has a right adjoint $r$ and $\phi$ is invertible then $\bar{u}$ has a right adjoint $\bar{r}$
with $x_s\bar{r}\cong rx_t$. 
\end{corollary}
 
 \section{Liftings}\label{liftings}
  
 We work in a bicategory $\mathfrak{C}$. We use the notation
\begin{equation}\label{rightlift}
 \begin{aligned}
   \xymatrix{ & & \CA \ar[dd]^{S} \\
     \CK  \ar[rru]^{\mathrm{rif}(S,B)} \ar[rrd]_{B} & & \dtwocell[0.4]{ll}{\varepsilon^S_B} \\
     & & \CB}
 \end{aligned}
\end{equation}
to depict a {\em right lifting} $\mathrm{rif}(S,B)$ (see \cite{12}) of the 1-morphism $B$ through the 1-morphism $S$.
The defining property is that pasting with \eqref{rightlift} gives a bijection 
$$\mathfrak{C}(\CK,\CA)(H, \mathrm{rif}(S,B))\cong \mathfrak{C}(\CK,\CB)(SH, B) \ .$$
(This concept is a dual of Kan extension.)

The right lifting is said to be {\em respected} by the 1-morphism $\CD\xra{K}\CK$ when $\varepsilon^S_B\cdot K$
exhibits $\mathrm{rif}(S,B)\cdot K$ as a right lifting of $B\cdot K$ through $S$. 
If $S$ has a right adjoint $T$ then $\mathrm{rif}(S,B) \cong T\cdot B$ and so is respected by all 1-morphisms $K$.
On the other hand, if $\mathrm{rif}(S,1_{\CB})$ exists and is respected by $S$ then
$S\dashv \mathrm{rif}(S,1_{\CB})$.     

If $B\xra{\alpha} C$ is a 2-morphism
and $B$ and $C$ have right liftings through $S$,
we write $\mathrm{rif}(S,B)\xra{\mathrm{rif}(S,\alpha)} \mathrm{rif}(S,C)$ for the 2-morphism 
defined by 
\begin{eqnarray*}
\left(S\mathrm{rif}(S,B)\xra{S\mathrm{rif}(S,\alpha)} S\mathrm{rif}(S,C)\xra{\varepsilon_C^S}C\right) =
\left(S\mathrm{rif}(S,B)\xra{\varepsilon_B^S} B\xra{\alpha}C\right) \ . 
\end{eqnarray*}

As with all cartesian morphisms, we have this simple property.  

\begin{proposition}\label{rifcancel} Suppose $\CA \xra{S} \CB \xra{U} \CC \xla{C}\CK$ are 1-morphisms
such that $\mathrm{rif}(U,C) : \CK\to \CB$ exists. Then
\begin{eqnarray*}
\mathrm{rif}(S,\mathrm{rif}(U,C)) \cong \mathrm{rif}(US,C)
\end{eqnarray*}
in the sense that one side exists if and only if the other does and $\varepsilon^{US}_B$ is the pasted
composite of $\varepsilon^{S}_{\mathrm{rif}(U,B)}$ and $\varepsilon^{U}_B$.  
\end{proposition}

\begin{proposition}\label{rifwithM} Suppose $\CA \xra{S}  \CB \xra{M} \CB \xla{B} \CK$ are 1-morphisms.
Suppose $1_{\CB}\xra{\eta} M$ is a 2-morphism with $\eta B$ a regular 
monomorphism in $\mathfrak{C}(\CK,\CB)$ and 
$\eta C$ a monomorphism in $\mathfrak{C}(\CK,\CB)$ for all $C\in \mathfrak{C}(\CK,\CB)$.
\begin{itemize}
\item[(i)]  For all $B\in \mathfrak{C}(\CK,\CB)$,
\begin{eqnarray}\label{eee}
\xymatrix @R-3mm {
B \ar@<0.0ex>[rr]^{\eta B} && MB \ar@<1.5ex>[rr]^{\eta MB}  \ar@<-1.5ex>[rr]_{M\eta B}  && M^2B 
}\end{eqnarray}
is an equalizer in $\mathfrak{C}(\CK,\CB)$.
\item[(ii)] Suppose right liftings $\mathrm{rif}(S,MB)$ and $\mathrm{rif}(S,M^2B)$ exist and are respected
by $\CD\xra{K}\CK$. 
The existence of a right lifting $\mathrm{rif}(S,B)$ respected by $\CD\xra{K}\CK$ is equivalent to the existence of an equalizer
\begin{eqnarray}\label{rifequalizer}
\xymatrix @R-3mm {
E \ar@<0.0ex>[rr]^{\kappa \ \phantom{aaa}} && \mathrm{rif}(S,MB) \ar@<1.5ex>[rr]^{\mathrm{rif}(S,\eta MB)}  \ar@<-1.5ex>[rr]_{\mathrm{rif}(S,M^2\eta B)}  && \mathrm{rif}(S,M^2B)  
}\end{eqnarray}
preserved by $\mathfrak{C}(\CK,\CA)\xra{\mathfrak{C}(K,1_{\CA})} \mathfrak{C}(\CD,\CA)$. 
\item[(iii)] In the situation of (ii), there is an isomorphism $\omega : E \cong \mathrm{rif}(S,B)$ whose composite with $\mathrm{rif}(S,\eta B)$ is $\kappa$.
\end{itemize}
\end{proposition} 
\begin{proof}
For (i), we know that $\eta B$ is the equalizer of some pair $\alpha, \beta : MB\to C$.
Take $\phi : D\to MB$ such that $\eta MB\cdot \phi=M\eta B\cdot \phi$. Then
\begin{eqnarray*}
\eta C\cdot \alpha \cdot \phi = M\alpha \cdot \eta MB\cdot \phi = M\alpha \cdot M\eta B\cdot \phi = M(\alpha \cdot \eta B)\cdot \phi 
\end{eqnarray*}
and similarly $\eta C\cdot \beta \cdot \phi = M(\beta \cdot \eta B)\cdot \phi$. 
So $\eta C\cdot \alpha \cdot \phi = \eta C\cdot \beta \cdot \phi$. Since $\eta C$ is a monomorphism,
we have $\alpha \cdot \phi = \beta \cdot \phi$ so that $\phi = \eta B \cdot \psi$ for some $\psi$ which
is unique because $\eta B$ is also a monomorphism. 

For (ii), using (i), we have
\begin{eqnarray}\label{equalizer=rif}
\lefteqn{\mathfrak{C}(\CK,\CA)(A,E) } \nonumber \\
& \cong & \{ A\xra{\sigma}\mathrm{rif}(S,MB) : \ \mathrm{rif}(S,\eta_{MB})\cdot\sigma = \mathrm{rif}(S,M\eta_{B})\cdot\sigma \} \nonumber   \\
& \cong & \{ SA\xra{\tau}MB : \ \eta_{MB}\cdot\tau= M\eta_{B}\cdot\tau  \}  \\
& \cong & \mathfrak{C}(\CK,\CB)(SA,B) \nonumber \\
& \cong & \mathfrak{C}(\CK,\CA)(A,\mathrm{rif}(S,B)) \nonumber \ .
\end{eqnarray}
naturally in $A\in \mathfrak{C}(\CK,\CA)$.

For (iii), by Yoneda, there is an isomorphism $\omega : E \cong \mathrm{rif}(S,B)$ inducing the composite
isomorphism \eqref{equalizer=rif}; this gives commutativity of the square
\begin{equation*}
\xymatrix{
SE \ar[rr]^-{S\kappa} \ar[d]_-{S\omega} && S\mathrm{rif}(S,MB) \ar[d]^-{\varepsilon_{MB}^S} \\
S\mathrm{rif}(S,B) \ar[r]_-{\varepsilon_{B}^S} & B \ar[r]_-{\eta B} & MB \ . } 
\end{equation*}
By applying the bijection 
$$\mathfrak{C}(\CK,\CB)(SE,MB)\cong \mathfrak{C}(\CK,\CA)(E,\mathrm{rif}(S,MB)) $$
we obtain the result stated. 
\end{proof}

\begin{proposition} Suppose $\CA \xra{S} \CB \xra{U} \CC$ and $\CK\xra{B}\CB$ are 1-morphisms.
Suppose $M : = \mathrm{rif}(U,U)$ exists and is respected by all 1-morphisms $\CK\to \CB$.
Suppose the 2-morphism $1_{\CB}\xra{\eta}M$, defined by $\varepsilon_U^U\cdot U\eta = 1_U$,  
is such that $\eta B$ is a regular monomorphism in $\mathfrak{C}(\CK,\CB)$ and 
$\eta C$ is a monomorphism in $\mathfrak{C}(\CK,\CB)$ for all $C\in \mathfrak{C}(\CK,\CB)$.
Suppose both $Q_B : = \mathrm{rif}(US,UB)$ and $Q_{MB} : = \mathrm{rif}(US,UMB)$ exist
and are respected by $\CD\xra{K}\CK$. 
Then $\mathrm{rif}(S,B)$ has the same universal property as the equalizer of the pair of
2-morphisms from $Q_B$ to $Q_{MB}$ in $\mathfrak{C}(\CK,\CA)$ corresponding 
(via the universal property of $Q_{MB}$) to the two paths in \eqref{eqrifSB}.
Moreover, $\mathrm{rif}(S,B)$ is respected by $K$ if and only if the equalizer is preserved
by $\mathfrak{C}(K,1_{\CA})$.   
\begin{equation}\label{eqrifSB}
\begin{aligned}
\xymatrix{
USQ_B \ar[d]_-{U\eta SQ_B} \ar[r]^-{\varepsilon^{US}_{UB}} & UB \ar[rr]^-{U\eta B}  & & UMB   \\
UMSQ_B \ar[r]_-{\cong}   &  U\mathrm{rif}(U,USQ_B)  \ar[rr]_-{U\mathrm{rif}(U,\varepsilon^{US}_{UB})} & & U\mathrm{rif}(U,UB) \ar[u]_-{\cong} }
\end{aligned}
\end{equation}
\end{proposition}
\begin{proof} By Proposition~\ref{rifcancel}, $\mathrm{rif}(S,MB)\cong\mathrm{rif}(S,\mathrm{rif}(U,UB)) \cong \mathrm{rif}(US,UB) = Q_B$ and $\mathrm{rif}(S,M^2B)\cong\mathrm{rif}(S,\mathrm{rif}(U,UMB)) \cong \mathrm{rif}(US,UMB) = Q_{MB}$; so they exist by assumption and Proposition~\ref{rifwithM}
applies. Furthermore, the parallel pair in \eqref{rifequalizer} transports across these isomorphisms to
the pair corresponding to the two paths in \eqref{eqrifSB}. 
\end{proof}

\begin{corollary}[Dubuc Adjoint Triangle Theorem \cite{DubucAT}] Suppose the 1-morphism $\CB \xra{U} \CC$ has a right adjoint $R$ with the unit $1_{\CB} \xra{\eta} RU$ a regular monomorphism preserved by the functors
$\mathfrak{C}(B,1_{\CB})$ for all 1-morphisms $B$ with target $\CB$. A 1-morphism $\CA \xra{S}\CB$ has a right adjoint if and only if the composite $US$ has a right adjoint $Q$ and the coreflexive pair of 2-morphisms
\begin{equation}\label{pairtoequalize}
\begin{aligned}
\xymatrix{
& QUSQU \ar[rr]^-{QU\eta {SQU}} & & QURUSQU \ar[rd]^-{QUR\alpha U} & \\
QU \ar[ru]^-{\beta QU} \ar[rrrr]_-{QU\eta}  & & & & QURU }
\end{aligned}
\end{equation}
admits an equalizer preserved by the functors $\mathfrak{C}(B,1_{\CB})$ for all 
1-morphisms $B$ with target $\CB$, where $\alpha$, $\beta$ are the counit and unit for $US\dashv Q$. 
\end{corollary}

\end{document}